\documentclass[12pt]{amsart}

\input{epsf}
\usepackage[english]{babel}
\usepackage[dvips]{graphicx}
\usepackage{amsmath,amssymb,amsfonts,latexsym}
\usepackage{graphics,psfig}

\newtheorem{theorem}{Theorem}[section]

\newtheorem{corollary}[theorem]{Corollary}
\newtheorem{lemma}[theorem]{Lemma}

\newtheorem{proposition}[theorem]{Proposition}

\newtheorem{example}[theorem]{Example}

\newtheorem{remark}[theorem]{Remark}

\def\A{\mathcal{A}}

\title{\textsc{On unavoidable sets of word patterns}}

\begin{document}

\maketitle
\thispagestyle{empty}
\begin{center}
Alexander Burstein\\
Department of Mathematics \\ Iowa State University\\
Ames, IA 50011-2064, USA\\
\texttt{burstein@math.iastate.edu}\\[1.3ex]
Sergey Kitaev\\
Department of Mathematics \\ University of Kentucky\\
Lexington, KY 40506-0027, USA\\
\texttt{kitaev@ms.uky.edu}\\[1.3ex]
\end{center}

\begin{abstract}

We introduce the notion of unavoidable (complete) sets of word
patterns, which is a refinement for that of words, and study
certain numerical characteristics for unavoidable sets of
patterns. In some cases we employ the graph of pattern overlaps
introduced in this paper, which is a subgraph of the de Bruijn
graph and which we prove to be Hamiltonian. In other cases we
reduce a problem under consideration to known facts on unavoidable
sets of words. We also give a relation between our problem and
intensively studied universal cycles, and prove there exists a
universal cycle for word patterns of any length over any alphabet.

\bigskip


\bigskip

\end{abstract}


\section{Introduction} \label{section_1}

When defining or characterizing sets of objects in discrete
mathematics, ``languages of prohibitions'' are often used to
define a class of objects by listing the prohibited subobjects,
i.e. subobjects that are not allowed to be contained in the
objects of the class. The notion of a subobject is defined in
different ways depending on the objects under consideration: a
subword (a block or segment) for fragmentarily restricted
languages, a subgraph for families of graphs, a subshape for
two-dimensional shapes (e.g. a submatrix for matrices) and so on.

We collect all prohibited objects into a set that we call a set of
prohibited objects, or simply a \emph{set of prohibitions}. The
idea of \emph{unavoidable} (or \emph{complete}\footnote{The word
``complete'' appears in e.g. ~\cite{evdok0}--\cite{evdkra}, but
the word ``\emph{unavoidable}'' is of common use in contemporary
literature (e.g. see~\cite[Chapter 3]{loth2}, \cite{sakhig}), so
we decided to use the latest terminology in this paper.}) set is
as follows: if there exists a restriction on the size of an
object, in other words, if large enough objects must contain
prohibited subobjects, then the set of prohibitions is
unavoidable.

In this paper, we are interested in unavoidable sets of \emph{word
patterns}, or just \emph{patterns} (see Section~\ref{section_3}
for definitions). These patterns are an extension of the
\emph{permutation patterns} studied extensively for the last
twenty years (see~\cite{kitman} for a survey on the corresponding
problems). Our unavoidable sets of patterns are refinements for
those of words. Questions on unavoidability of sets of words
appear, for instance, in algebra (sequences without repetitions),
coding theory (chain codes), number theory (arithmetic
progressions in partitions of the set of natural numbers),
dynamical systems (motions of an object in a space with certain
restrictions).

There is a number of numerical characteristics that are valuable
for unavoidability criteria and the recognition algorithms based
on them. Three such characteristics, namely $M_w(n)$, $L_w(n)$ and
$C_w(n)$ (for definitions see Section~\ref{section_2}), are
considered in~\cite{evdok1}. We consider the similar
characteristics $M_p(n,m)$, $L_p(n,m)$ and $C_p(n,m)$ for the case
of prohibited patterns (for definitions see
Section~\ref{section_3}), where $m$ is the number of letters in
the corresponding alphabet (we do not use this parameter for the
functions $M_w(n)$, $L_w(n)$ and $C_w(n)$ to be consistent
with~\cite{evdok1}). Moreover, in Subsubsection~\ref{cycle} we
discuss how finding a lower bound for $C_p(n,m)$ is related to the
so-called \emph{universal cycles for combinatorial structures}
that have been studied intensively (e.g. see~\cite{chung,hur} and
references therein). To get the lower bound, we prove that the
\emph{graph of pattern overlaps} (see definition in
Section~\ref{section_3}) is Hamiltonian, and derive as a corollary
that there exists a universal cycle for word patterns of any
length over any alphabet (see Corollary~\ref{corcor}).

We remark that when considering patterns, the underlying alphabet
must be ordered, as opposed to the objects considered
in~\cite{evdok1}.

The paper is organized as follows. In Section~\ref{section_2} we
review the main results on unavoidable sets of words
in~\cite{evdok1,evdok2}. The motivation for a relatively detailed
review of these papers is the fact that they are available only in
Russian (as far as we know), which caused, in particular, the
rediscovery of some of those results in~\cite{sakhig}. Besides,
the results obtained in~\cite{evdok1,evdok2} are of great interest
in general and very useful in this paper in particular. In
Section~\ref{section_3}, we define the notion of a \emph{pattern},
an \emph{$n$-pattern word}, and study unavoidable sets of
patterns.

\section{Unavoidable sets of words} \label{section_2}

Let $\A =\{ a_1, \ldots ,a_n\}$ be an alphabet of $n$ letters. A
\emph{word} over the alphabet $\A$ is a finite sequence of letters
of the alphabet. Any $i$ consecutive letters of a word $X$
generate a \emph{subword} of length $i$. The set ${\A}^*$ is the
set of all words over the alphabet ${\A}$, and ${\A}^n$ is the set
of all words over $\A$ of length $n$. Let $S \subseteq {\A}^*$ be
a set of prohibited words or a set of prohibitions. A word that
does not contain any words from $S$ as its subwords is said to be
\emph{free} from $S$ or \emph{$S$-free}. The set of all $S$-free
words is denoted by $\widehat{S}$.

If there exists a natural number $k$ such that the length of any
word in $\widehat{S}$ is less than $k$, then $S$ is called an
unavoidable set. This is straightforward to see that $S$ is
unavoidable if and only if $\widehat{S}$ has finitely many of
elements. Thus, for any unavoidable set $S$ we can define the
function
\[
L_w(\widehat{S})=\max_{X\in\widehat{S}}{\ell(X)},
\]
where $\ell(X)$ is the length of a word $X$.

The basic problem in considering of sets of prohibitions is
whether or not a given set $S$ of prohibitions is unavoidable.
Other possible questions are: given an unavoidable $S$ find or
estimate $L_w(\widehat{S})$; construct an $S$-free word of length
$L_w(\widehat{S})$; find the number of elements in $\widehat{S}$.
If $S$ is avoidable then some possible questions are: find an
infinite $S$-free sequence; describe all such sequences; find the
cardinality of the set of these sequences; find the cardinality of
the set of finite $S$-avoiding sequences of a given length.

Let $S$ be a finite set of words over an alphabet $\A$, and let
$n$ be the maximal length of a word in $S$. If a word $X$ is a
subword of a word $Y$ then we say that $Y$ is a \emph{superword}
for $X$. Suppose now that a word $X\in S$ and $\ell(X)<n$. Remove
$X$ from $S$ and adjoin to $S$ all superwords for $X$ of length
$n$. If this procedure is performed for any such $X$, and all
resulting repetitions are removed, we will get a set $S^{'}$ of
distinct words of length $n$.

\begin{proposition}\label{same_lengths}
{\rm(}\cite[Proposition 1]{evdok1}{\rm)} $S$ is unavoidable iff
$S^{'}$ is unavoidable.
\end{proposition}

Thus, sets of prohibitions $S\subseteq {\A}^n$ are of special
interest, and for the most part, our considerations in this paper
are related to these sets. More precisely, we will consider the
functions
\[
M_w(n)=\min |S| \quad \mbox{and} \quad L_w(n) = \max
L_w(\widehat{S}),
\]
where the extremum is taken with respect to all unavoidable
$S\subseteq {\A}^n$. These functions are examples of numerical
characteristics that describe the bound between avoidable and
unavoidable sets of prohibitions. To give an instance of such a
bound, we consider the following example.

\begin{example}\rm(\cite[Examples 1,2]{evdok1}).
Consider $\A=\{0,1\}$ and the sets of prohibitions
\[
\begin{split}
S_1 = \{ 000, 001, 101\underline{1}, 0101, 1111 \},\\
S_2 = \{ 000, 001, 101\underline{0}, 0101, 1111 \}.
\end{split}
\]
Thus $S_1$ and $S_2$ differ only in one underlined letter. One can
see that $S_1$ is unavoidable, and $L_w(\widehat{S_1})=8$. On the
other hand, $S_2$ is avoidable. Indeed,
\[
\underbrace{011}\underbrace{011}\ldots \quad \mbox{and} \quad
\underbrace{0111}\underbrace{0111}\ldots
\]
are $S_2$-free, and
\[
\underbrace{011}\underbrace{0111} \quad \mbox{and} \quad
\underbrace{0111}\underbrace{011}
\]
are $S_2$-free. Hence, substituting $0\mapsto 011$ and $1\mapsto
0111$ in any sequence over $\A$, we get an $S_2$-free sequence.
Hence, the cardinality of $\widehat{S_2}$ is the continuum.
\end{example}

\begin{center}
\begin{figure}[h]
\epsfxsize=200.0pt
\epsffile{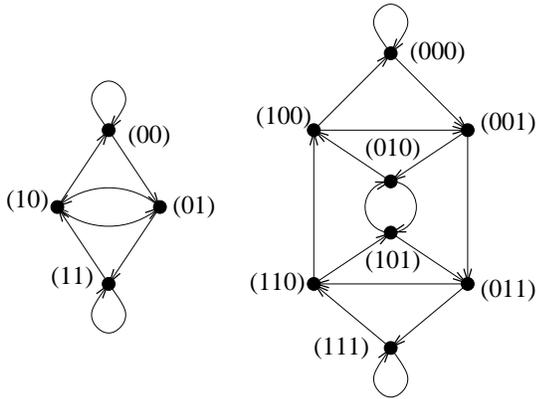}
\caption{The de Bruijn graphs for the alphabet $\A=\{0,1\}$ and $n=2,3$.}
\label{fig1}
\end{figure}
\end{center}

In what follows, we will need the following graph. A \emph{de
Bruijn graph} is a directed graph $\vec{G}_n=\vec{G}_n(V,E)$,
where the set of vertices $V$ is the set of all words in $\A^n$,
and there is an arc from $u\in\A^n$ to $v\in\A^n$ if and only if
\[
u=aw \text{ and } v=wb \quad \text{for some } w\in\A^{n-1} \text{
and } a,b\in\A.
\]
Figure~\ref{fig1} shows the de Bruijn graphs for a 2-letter
alphabet and $n=2,3$.

The de Bruijn graphs were first introduced (for the alphabet
$\A=\{0,1\}$) by de Bruijn in 1944 for finding the number of code
cycles. However, these graphs proved to be a useful tool for
various problems related to combinatorics on words (e.g.
see~\cite{evdok1, evdok2, golomb}). It is known that the graph
$\vec{G}_n$ can be defined recursively as
$\vec{G}_n=L(\vec{G}_{n-1})$, where $L$ indicates the operation of
taking the line graph.

\begin{center}
\begin{figure}[h]
\epsfxsize=80.0pt \epsffile{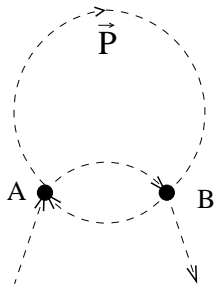} \caption{The arc $\vec{BA}$
is a chord for the path $\vec{P}$, but $\vec{AB}$ is not.}
\label{fig2}
\end{figure}
\end{center}

A \emph{chord} of a directed simple path $\vec{P}$ in $\vec{G}_n$
is an arc that does not belong to $\vec{P}$ but connects two of
its vertices in a such way that there is a circuit generated by
this arch and the part of the path between the ends of the arc.
For instance, on Figure~\ref{fig2} the arc $\vec{BA}$ is a chord
for the path $\vec{P}$, whereas $\vec{AB}$ is not.

Let $C_w(n)$ denote the greatest length (the number of vertices)
of a simple path in $\vec{G}_n$ that does not have chords and does
not go through any vertex that has a loop. The following theorem
was proved by considering the de Bruijn graph.

\begin{theorem}\label{thpath}{\rm(\cite[Theorem 1]{evdok1})}
$L_w(n) = C_w(n) + n - 1 = |\A|^{n-1} + n - 2.$
\end{theorem}

The following theorem was proved using the cyclic structure of the
de Bruijn graph (the main result of~\cite{golomb}) as well as the
number of conjugacy classes of words with respect to a cyclic
shift.

\begin{theorem}\label{theorphi} {\rm(\cite[Theorem 2]{evdok1})}
\[
M_w(n) = \frac{1}{n}\sum_{d|n} \varphi (n/d) |\A|^d,
\]
where $\varphi(n)$ is the number of integers in
$\{1,2,\ldots,n-1\}$ relatively prime to $n$ (Euler's
$\varphi$-function).
\end{theorem}

Since any set of prohibitions $S$ with $|S|< M_w(n)$ is avoidable,
it is helpful to have a table for $M_w(n)$. For $|\A| = 2$ and $2
\le n \le 10$, see Table~\ref{table1}.

\begin{table}[ht]
\begin{center}
\begin{tabular}{c|c|c|c|c|c|c|c|c|c}
$n$ & 2 & 3 & 4 & 5 & 6 & 7 & 8 & 9 & 10 \\
\hline
$M_w(n)$ & 3 & 4 & 6 & 8 & 14 & 20 & 36 & 60 & 108 \\
\end{tabular}
\caption{The function $M_w(n)$ for $2\le n \le 10$ and a 2-letter
alphabet.} \label{table1}
\end{center}
\end{table}

In particular, any set of binary words of length 9 that has less
than 60 words is avoidable. Also, it is obvious that $M_w(n)\sim
|\A|^n/n$, when $n \to \infty$. The last observation allows us to
prove the following statement.

\begin{proposition}{\rm(\cite[Proposition 1]{evdok2})}
There exist at least $2^{|\A|^n(1-{\varepsilon}_n)}$ unavoidable
sets $S \subseteq {\A}^n$. Here ${\varepsilon}_n \to 0$ when $ n
\to \infty$.
\end{proposition}

\section{Unavoidable sets of patterns}\label{section_3}

The alphabets considered in this section must be totally ordered,
and without loss of generality they coincide with
$[m]=\{1,2,\ldots,m\}$ for an appropriate $m$.

We refer to~\cite{kitman} for a general survey of various pattern
problems. However, in this paper we are concerned only with
\emph{word patterns} studied for the first time in \cite{burman}.
More precisely, we consider the word patterns \emph{without
internal dashes} (see~\cite{kitman}). For this paper, we can
define a pattern to be a subword (of a word) that contains each of
the letters $1,2,\ldots,k$ at least once for some $k$, and no
other letters. For instance, the word 2613235 contains an
occurrence of the pattern $1323$, but its subword 2613 is not a
pattern. By analogy with Section~\ref{section_2}, if a word does
not contain a pattern $p$, it is \emph{free} from $p$ or
\emph{$p$-free}. However, the crucial difference between this
section and Section~\ref{section_2} is that instead of considering
\emph{words} free from a pattern $p$, we consider the objects that
we call the \emph{$n$-pattern words}. An $n$-pattern word is a
word in which each subword of length $n$ is a pattern. Thus,
constructing $n$-pattern words, we can restrict ourselves to
alphabets having at most $n$-letters. Indeed, an occurrence of a
letter $m>n$ in a subword $A$ of length $n$ of an $n$-pattern word
$W$ contradicts the fact that $A$ must be a pattern ($A$ must
contain each of the letters $1,2,\ldots,m$).

By analogy with Section~\ref{section_2}, when dealing with sets of
prohibited words, we can consider sets of prohibited patterns, or
simply sets of prohibitions, when it is clear which prohibitions
we mean. We can also define the notion of an unavoidable set here
in the same way. However, in considering prohibited patterns and
$n$-pattern words, we assume that all prohibitions are of length
$n$. Hence, for patterns, we can define the functions $L_p(n,m)$
and $M_p(n,m)$ similarly to $L_w(n)$ and $M_w(n)$ (recall that $m$
is the number of letters in the alphabet). As in
Section~\ref{section_2}, the basic problem is whether or not a
given set $S_p$ of prohibitions is unavoidable, and $L_p(n,m)$ and
$M_p(n,m)$ are important numerical characteristics to study.

\subsection{The function $M_p(n,m)$}

Recall that the M\"obius function is defined by
\[
\mu (n) =
\begin{cases}
0, & \text{if $n$ has one or more repeated prime factors},\\
1, & \text{if $n=1$}, \\
(-1)^k, & \text{if $n$ is a product of $k$ distinct primes},
\end{cases}
\]
so $\mu(n)\ne 0$ indicates that $n$ is square-free.

The purpose of this subsection is to prove the following theorem.

\begin{theorem}\label{main1}
For $n$-pattern words over $[m]$, we have
\[
M_p(n,m)=\sum_{i|n}\sum_{j=0}^{\min(i,m)-1}
(-1)^j\binom{\min(i,m)-1}{j}\frac{1}{i}
\sum_{d|i}\mu(d)(\min(i,m)-j)^{\frac{i}{d}},
\]
where $M_p(n,m)=\min|S_p|$, and the minimum is taken over all
unavoidable sets $S_p$ of patterns of length $n$ over the alphabet
$[m]$.
\end{theorem}

One can compare this result with that of Theorem~\ref{theorphi}.

\begin{remark} \rm In Theorem~\ref{main1}, we can assume that $n\ge m$,
since if $n<m$ we can only use the first $n$ letters in $[m]$ to
construct $n$-pattern words, which reduces to the case $n=m$.
\end{remark}

\begin{remark} \rm For $n=m$, we have $\min(i,m)=i$ in the formula of
Theorem~\ref{main1}.
\end{remark}

\begin{center}
\begin{figure}[h]
\epsfxsize=200.0pt \epsffile{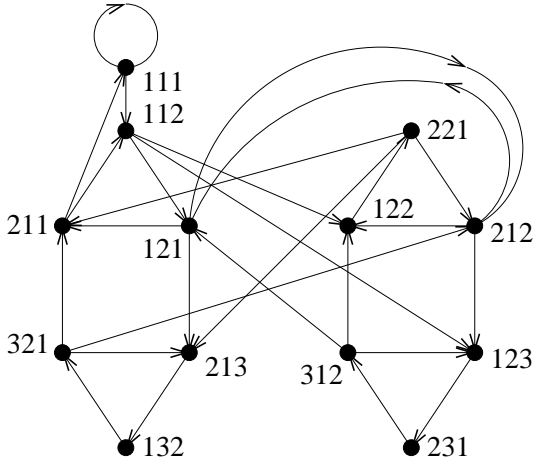} \caption{The graph
of pattern overlaps for $\A=\{1,2,3\}$ and $n=3$.} \label{fig4}
\end{figure}
\end{center}

To prove Theorem~\ref{main1}, we introduce the \emph{graph of
pattern overlaps} $\vec{P}_n=\vec{P}_n(V,E)$, which is a subgraph
of the de Bruijn graph $\vec{G}_n$, where the set of vertices $V$
contains all $n$-letter patterns over the underlying alphabet
$\A$, and the set of arcs $E$ consists of all the arcs of
$\vec{G}_n$ between vertices corresponding to the patterns. In
Figure~\ref{fig4}, we can see the graph of pattern overlaps in the
case of a 3-letter alphabet and $n=3$ (we omit parentheses around
the triples on the graph to indicate that we are dealing with
$\vec{P}_3$, not $\vec{G}_3$).

Let $T_p(n,m)$ denote the number of conjugacy classes of patterns
of length $n$ over the alphabet $[m]$ with respect to a cyclic
shift. For instance, there are 5 conjugacy classes on
Figure~\ref{fig4}. They are $\{111\}$, $\{112, 121, 211\}$,
$\{221, 212, 122\}$, $\{321, 213, 132\}$ and $\{312, 123, 231\}$.
Thus, $T_p(3,3)=5$.

\begin{lemma} $M_p(n,m)=T_p(n,m)$ \label{lemma1}
\end{lemma}

\begin{proof}
To prove the lemma, we follow the proof of Theorem~\ref{theorphi}
in~\cite{evdok1}.

Suppose $S_p$ is an unavoidable set of patterns of length $n$ and
$X$ is an arbitrary $n$-pattern word of length $n$ ($X$ is a
pattern) over $[m]$. We form the sequence
\[
X^{\infty}=X X X \ldots,
\]
by repeating the word $X$ periodically. Since $S_p$ is
unavoidable, $X^{\infty}$ contains a prohibited pattern $p\in
S_p$. From the construction of the sequence, $p$ is either $X$ or
a cyclic shift of $X$. Thus $S_p$ contains a pattern from each
conjugacy class of patterns of length $n$ over $[m]$ with respect
to a cyclic shift. Thus, $|S_p|\ge T_p(n,m)$, and since $S_p$ is
an arbitrary set, we have
\[
M_p(n,m)\ge T_p(n,m).
\]

To prove that $T_p(n,m)$ is an upper bound, we need to find an
unavoidable set of cardinality $T_p(n,m)$. We consider the graph
$\vec{P}_n$ whose vertices correspond to the words over $[m]$. If
$V'\subset V(\vec{P}_n)$ and each circuit of $\vec{P}_n$ contains
a vertex in $V'$ then we say that $V'$ \emph{cuts} all circuits of
$\vec{P}_n$. By deleting all such $V'$ with all incident arcs from
$\vec{P}_n$, we get an acyclic graph on the vertex set
$V\backslash V'$. The set of the patterns in $[m]^n$ corresponding
to the vertices in $V'$ is unavoidable. Indeed, if not, a sequence
free from $V'$ determines a self-intersecting walk in $\vec{P}_n$
and thus generates a circuit on the vertex set $V\backslash V'$,
which is impossible.

Golomb~\cite{golomb} found a set of vertices $V_c$ that cuts all
circuits of the de Bruijn graph $\vec{G}_n$ with $|V_c|$ equal to
the number of conjugacy classes of the words. Thus $V_c$ cuts all
circuits in $\vec{G}_n$ and has one vertex in each conjugacy
class. Since $\vec{P}_n$ is a subgraph of $\vec{G}_n$, $\vec{P}_n$
will have no circuit after removing the vertices in $V_c$. The set
of vertices in $V_c$ that belong to $\vec{P}_n$ corresponds to an
unavoidable set, and thus
\[
M_p(n,m)\le T_p(n,m).
\]
This proves the lemma.
\end{proof}


\begin{lemma}\label{lemma2}
\[
T_p(n,m)=\displaystyle\sum_{i|n} \sum_{j=0}^{\min(i,m)-1}
(-1)^j\binom{\min(i,m)-1}{j}\frac{1}{i}
\sum_{d|i}\mu(d)(\min(i,m)-j)^{\frac{i}{d}}.
\]
\end{lemma}

\begin{proof}
Recall that a word $x\in\A^*$, where $\A$ is any (ordered or
unordered) alphabet, is called \emph{primitive} if it is not a
power of another word. Thus $x\ne\emptyset$ is primitive if
$x=y^e$ only for $e=1$. For instance, the words 121, 1221, 12121
are primitive, whereas the word 121212 is not. It is easy to show
that each nonempty word is a power of a unique primitive word.
Thus, $x=r^e$ for a unique primitive word $r$. The number $e$ is
called the \emph{exponent} of $x$. It is also easy to see that all
words, and hence all patterns, in the same conjugacy class have
the same exponent. Moreover, if $x_1=r_1^e$ and $x_2=r_2^e$ and
$|x_1|=|x_2|$, then $x_1$ is conjugate to $x_2$ iff $r_1$ is
conjugate to $r_2$. We define the notion of a \emph{primitive
pattern} in the same way as for words. Clearly, all properties of
primitive words hold for primitive patterns as well.

So, in order to find $T_p(n,m)$, we need to find the number of
conjugacy classes of primitive patterns of length $i$ over the
alphabet $[m]$, where $i|n$, and then take a sum of these numbers.
However, for a given $i$, we cannot use directly the well known
formula for the number of conjugacy classes of primitive words
over $\min(i,m)$-letter alphabet (a primitive word of length $i$
can have at most $i$ distinct letters, since we are dealing with
patterns), given by
\[
\frac{1}{i}\sum_{d|i}\mu(d)(\min(i,m))^{\frac{i}{d}}.
\]
Indeed, this formula counts, among others, primitive words which
are not primitive patterns (when some letter $j$, $2\le j\le
\min(i,m)-1$, occurs in a primitive pattern whereas $j-1$ does
not). So, we need to use the standard inclusion-exclusion method
(the sieve formula) to handle this situation. We define the
property $A_j$ to be ``the letter $j$ does not occur in a
primitive word''. Clearly we may restrict ourselves to the case
$j\le \min(i,m)-1$, since the absence of the largest letter,
namely $\min(i,m)$, is not a bad property when considering
patterns. Now we easily get the number of primitive patterns of
length $i$, which is given by
\[
\sum_{j=0}^{\min(i,m)-1}(-1)^j\binom{\min(i,m)-1}{j}
\frac{1}{i}\sum_{d|i}\mu(d)(\min(i,m)-j)^{\frac{i}{d}}.
\]
This proves the lemma.
\end{proof}

Now the truth of Theorem~\ref{main1} follows from Lemmas~\ref{lemma1}
and~\ref{lemma2}.

\subsection{The function $L_p(n,m)$}

Let $C_p(n,m)$ denote the greatest length (the number of vertices)
of a simple path in $\vec{P}_n$ that does not have chords (see the
definition in Section~\ref{section_2}) and does not pass through
any vertex incident with a loop. Using exactly the same
considerations as in the proof of Theorem~\ref{thpath}
(see~\cite{evdok1}), one can prove the following theorem.

\begin{theorem}\label{ththth}
$L_p(n,m) = C_p(n,m) + n - 1$.
\end{theorem}

Moreover, in the case $m=2$, the de Bruijn graph $\vec{G}_n$
almost coincides with the graph of pattern overlaps $\vec{P}_n$.
Indeed, the only difference between these graphs is the vertex
$(22\ldots 2)$ and all edges adjacent to that vertex ($22\ldots 2$
is the only binary non-pattern). However, the lemma
to~Theorem~\ref{thpath} (see~\cite{evdok1}) provides that in the
binary case $C_w(n)=2^{n-1}-1$, and since $C_w(n)$ is the maximal
length of a path that, in particular, does not pass through the
loop $(22\ldots 2)$, we have that in this case $C_w(n)=C_p(n,2)$.
Thus the following theorem is true:

\begin{theorem}\label{thth1}
$L_p(n,2)=2^{n-1}+n-2.$
\end{theorem}

However, in the case $m\ge 3$, the only useful information we can
extract from Theorem~\ref{thpath} is the following rough bound
\[
L_p(n,m)<m^{n-1}+n-2.
\]
So, according to Theorem~\ref{ththth} we need to find $C_p(n,m)$
in order to get $L_p(n,m)$. The purpose of the rest of the
subsection is to find an upper and a lower bound for $C_p(n,m)$
for $m\ge 3$.

\subsubsection{An upper bound for $C_p(n,m)$}
We only give a trivial upper bound. Clearly, in order to avoid
chords, each conjugacy class (with respect to shift) which has $i$
words can have no more than $i-1$ words in the path. Thus, we use
the formula for $T_p(m,n)$ with a correction, namely the factor of
$i-1$, which indicates that each primitive word of length $i$ is
responsible for a conjugacy class of $i$ elements, and we take
$i-1$ elements out of these $i$:
\[
C_p(n,m)\le \sum_{i|n}(i-1)
\sum_{j=0}^{\min(i,m)-1}(-1)^j\binom{\min(i,m)-1}{j}\frac{1}{i}
\sum_{d|i}\mu(d)(\min(i,m)-j)^{\frac{i}{d}}.
\]

\subsubsection{A lower bound for $C_p(n,m)$}\label{cycle}
We observe that the line graph $L(\vec{P}_{n-1})$ for the graph
$\vec{P}_{n-1}$ determines a subgraph of the graph $\vec{P}_n$. We
get that by using the general properties of the de Bruijn graph
(since $\vec{P}_n$ is its subgraph), as well as the fact that if
$x_1x_2\ldots x_{n-1}$ and $x_2x_3\ldots x_n$ are vertices in
$\vec{P}_{n-1}$, then the arc between them generates the vertex
$x_1x_2\ldots x_n$ in the line graph, and $x_1x_2\ldots x_n$ is a
pattern and thus belongs to $\vec{P}_n$. Moreover, from the
considerations in the proof of Theorem~\ref{thpath}
(see~\cite{evdok1}), it follows that a simple path in
$\vec{P}_{n-1}$ determines a simple path without chords in
$\vec{P}_n$ after removing the loop $11\ldots 1$.

So, in order to get a lower bound for $C_p(n,m)$, we need to
construct a simple path in $\vec{P}_{n-1}$ of as great a length as
possible (ideally a Hamiltonian path). In order to get a
Hamiltonian path or a path that is ``close'' to a Hamiltonian one,
we can try to use the methods and techniques similar to those used
in constructions of universal cycles for various combinatorial
structures such as words, permutations, partitions, and others
(e.g. see \cite{chung,hur}).

We briefly discuss the general notion of a universal cycle
(see~\cite{chung}).

Suppose we are given a family $\mathcal{F}_n$ of combinatorial
objects of ``rank $n$'' and let $m:=|\mathcal{F}_n|$ denote their
number. We assume that each $F\in\mathcal{F}_n$ is ``generated''
or specified by some sequence $x_1x_2\ldots x_n$, where
$x_i\in\mathcal{A}$ for some fixed alphabet~$\mathcal{A}$. We say
that $U=a_0a_1\ldots a_{m-1}$ is a universal cycle (or a
$U$-\emph{cycle}) for $\mathcal{F}_n$ if $a_{i+1}a_{i+2}\ldots
a_{i+n}$, $0\le i< m$, runs through each element of
$\mathcal{F}_n$ exactly once, where index addition is performed
modulo $n$.

In our case the combinatorial objects are patterns of length $n$,
and as in many other cases (e.g. \emph{de Bruijn cycles},
permutations, partitions), but not in all cases (e.g. $k$-subsets
of an $n$-set), it is possible to define a directed
\emph{transition graph}, namely the graph of pattern overlaps
$\vec{P}_n$, and reduce the problem of constructing a U-cycle to
constructing a Hamiltonian circuit for $\vec{P}_n$. Even though we
do not need a Hamiltonian circuit (since we are concerned with
paths of maximal length), but we can still try to use the same
techniques as in~\cite{chung,hur} and in references therein.

However, it turns out that the abovementioned techniques work only
for $m=2$, which we are not interested in since we have an
explicit result in this case (see Theorem~\ref{thth1}). The main
problem is that the graph of pattern overlaps is not balanced,
i.e. we have vertices where the indegree is not equal to the
outdegree. Also, $\vec{P}_n$ is not the line graph of
$\vec{P}_{n-1}$. However, it is possible to prove the following
statement.

\begin{theorem}\label{pattgraph}
The graph of pattern overlaps $\vec{P}_n$ contains a Hamiltonian
circuit.
\end{theorem}
\begin{proof}
We first observe that $\vec{P}_n$ is strongly connected. Indeed,
suppose we are given two vertices of $\vec{P}_n$, namely
$X=x_1x_2\ldots x_n$ and $Y=y_1y_2\ldots y_n$. If $I$ denotes the
vertex $11\ldots 1$, then we can find a path $\vec{P}_X$ from $X$
to $I$. Indeed, If $x_i$ is the largest letter in $X$, then we
consider the following path in $\vec{P}_n$:
\[
X=x_1x_2 \ldots x_n\to x_2x_3\ldots x_nx_1\to\cdots
\to x_ix_{i+1}\ldots x_{i-1}\to x_{i+1}\ldots x_{i-1}1=X'.
\]
Thus, in $X'$ we get 1 in place of the largest letter of $X$. We
observe that $X'$ is obviously a pattern. Clearly, we can continue
this path by replacing the largest letters, one by one, with 1's
until we arrive at $I$. On the other hand, it is easy to see that
the operation of changing a largest letter to 1 is invertible. For
instance, in order to find a path from $X'$ to $X$, we may do the
following sequence of steps:
\[
X'=x_{i+1}\ldots x_{i-1}1\to x_{i+2}\ldots 1x_{i+1} \to
\cdots\to 1x_{i+1}\ldots x_nx_1\ldots x_{i-1}\to$$
$$x_{i+1}\ldots x_nx_1\ldots x_{i-1}x_i\to
x_{i+2}\ldots x_ix_{i+1}\to \cdots \to x_1x_2\ldots x_n=X.
\]
Thus, we can find a path from $I$ to $Y$, which together with the
path $\vec{P}_X$, gives a path from $X$ to $Y$. Similarly, one can
get a path from $Y$ to $X$, which proves that $\vec{P}_n$ is
strongly connected.

The main property we use when proving $\vec{P}$ has a Hamiltonian
circuit is illustrated in Figure~\ref{fig5}A. It says that if
$C_1$ and $C_2$ are two circuits corresponding to different
conjugacy classes with respect to the shift, and there is an arc
from $C_1$ to $C_2$ then there is an arc from $C_2$ to $C_1$ and
vise versa. Moreover, in all cases but one (see discussion below),
we can choose these arcs as in the Figure~\ref{fig5}A, that is
once we leave $C_1$ at the vertex $xW$, we can come back, after
visiting $C_2$, at the vertex $Wx$, which is adjacent to $xW$ on
the circuit $C_1$. The notation $xW$ (resp. $Wx$) is used to
indicate a pattern of length $n$ with the first (resp. last)
letter $x$. The only exception when the picture differs from that
on Figure~\ref{fig5}A is the loop $11\ldots 1$, and there is only
circuit adjacent to it, namely the one generated by $11\ldots 12$.
In this case $xW$ coincides with $Wx$, which however does not
affect our considerations below.

\begin{center}
\begin{figure}[h]
\epsfxsize=370.0pt \epsffile{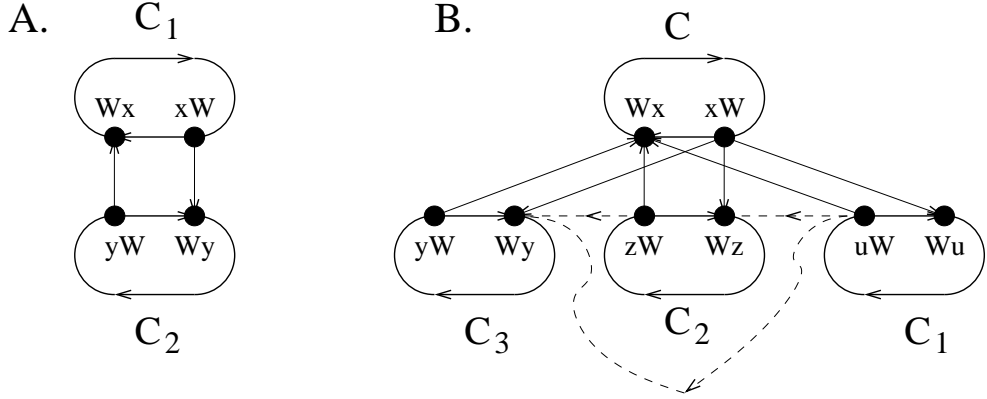} \caption{Circuits in
$\vec{P}_n$.} \label{fig5}
\end{figure}
\end{center}

\textbf{The basic idea:} We show the existence of a Hamiltonian
circuit iteratively, starting from any circuit corresponding to a
conjugacy class with respect to the shift, and on each following
iteration creating a new circuit that contains the previous one
and has more vertices since it covers additional circuits
corresponding to some conjugacy classes (by covering here we mean
containing all the vertices from a circuit in our big circuit).
Moreover, we construct the big circuit so that once it arrives at
a new circuit corresponding to a conjugacy class, it uses all the
vertices from that circuit before leaving. We keep doing that
using the fact that $\vec{P}_n$ is a disjoint union of the
circuits corresponding to the conjugacy classes, until we create a
Hamiltonian circuit.

Let $H_1$ be an arbitrary circuit corresponding to a conjugacy
class with respect to the shift. Now assume we made $i$ iterative
steps and obtained a circuit $H_i$. If $H_i$ covers all the
vertices of $\vec{P}_n$, then we are done. Otherwise, on iteration
$i+1$ we proceed as follows.

The fact that $\vec{P}_n$ is strongly connected ensures that there
is an arc from a circuit $C$ covered by $H_i$ to a circuit which
is not covered by $H_i$. Our strategy is to start from the vertex
where $H_i$ arrived at $C$, then go around $C$ following $H_i$
vertex by vertex, until we reach the vertex in which $H_i$ leaves
$C$, and at each step, checking if it is possible to extend $H_i$
according to the following considerations.

Assume we are in the vertex $xW$ in $C$. If there is only one arc
coming out of $C$, namely the arc to the vertex $Wx$ belonging to
$C$, then we cannot extend $H_i$ at this step, so we need to
consider the next vertex $Wx$ instead. Otherwise, there are $j\ge
1$ arcs that come out from $xW$ to $j$ different circuits
corresponding to some conjugacy classes (we denote the set of
these circuits by $B$). The case $j=1$ is shown on
Figure~\ref{fig5}A, if we assume $C=C_1$. In this case there are
two possibilities: either $C_2$ is covered by $H_i$ or not. In the
first case we cannot extend $H_i$, so we need to consider the
vertex $Wx$ belonging to $C$ to proceed further. In the second
case, we can extend $H_i$ by going to the vertex $Wy$, then
through the vertices belonging to $C_2$ until we reach $yW$, then
we come back to $C$ at the vertex $Wx$.

When $j>1$, either all circuits from $B$ are already covered by
$H_i$, or there is a number of circuits that are not covered by
$B$ $B$ (we denote the set of these circuits by $B_0$). In the
first cannot extend $H_i$ and we need to continue to proceed to
the vertex $Wx$. We claim that in the second case there is a path
starting from the vertex $xW$, going through all the vertices from
the circuits from $B_0$ and coming to the vertex $Wx$. We can
extend $H_i$ with this path. This claim is not hard to prove for
any $j$, for instance by induction. However, we only give our
proof in the case $j=3$ (see Figure~\ref{fig5}B) as it is easily
generalizable.

In Figure~\ref{fig5}B, $Wy$, $Wz$, and $Wu$ are representatives
from the circuits $C_3$, $C_2$ and $C_1$ respectively, which
belong to $B$. The key observation here is that any other circuit
in $B$ is as good as $C$, that is, e.g. we can go from $xW$ to any
of the vertices $yW$, $zW$ and $uW$, but we can also go from, say,
$uW$ to any of these vertices. If $B=B_0$, then we can start at
$xW$, go to $Wu$, go to $uW$ through $C_1$, then to $Wz$, then go
to $zW$ through $C_2$, to $Wy$, to $yW$ through $C_3$ and finally
come to $Wx$, in which case we succeeded to extend $H_i$. If $B\ne
B_0$, we use the same procedure simply skipping the circuits not
in $B_0$. E.g. if $C_2\notin B_0$, we change the path above by
going from $uW$ directly to $Wy$, again extending $H_i$.

Thus, we constructed the circuit $H_{i+1}$ that contains more
vertices than $H_i$ does. Since $\vec{P}_n$ has finitely many
vertices, $\vec{P}_n$ must contain a Hamiltonian circuit.
\end{proof}

\begin{remark} \rm
The proof of theorem~\ref{pattgraph} can be simplified, if we add
exactly one circuit corresponding to a conjugacy class at each
iteration. Indeed, in this case we do not need to consider the
sets $B$ and $B_0$ used in the proof, as well as the illustration
on Figure~\ref{fig5}B. Thus, once we find a circuit to add to the
big circuit, we can start a new iteration. However, we keep the
more complicated proof since it helps understand the structure of
the graph of pattern overlaps more deeply.
\end{remark}

\begin{remark} \rm
One can test how the algorithm of finding a Hamiltonian circuit in
$\vec{P}_n$ works in the case $n=3$ and $m=3$ on
Figure~\ref{fig4}.
\end{remark}

As an immediate corollary to Theorem~\ref{pattgraph} we have the
following:

\begin{corollary}\label{corcor}
For any $m$ and $n$, there exists a U-cycle for word patterns of
length $n$ over an $m$-letter alphabet.
\end{corollary}

The following proposition is easy to prove using elementary
combinatorics.

\begin{proposition}\label{propro} The number of different word
patterns of length $n$ on $m$ letters is
\[
\sum_{i=1}^{m}\sum_{a_1+\cdots+a_i=n \atop
a_1\ge1,\dots,a_i\ge1}\binom{n}{a_1,\ldots,a_i}.
\]
\end{proposition}

Now, using the discussion in the beginning of the subsubsection,
Theorem~\ref{pattgraph} and Proposition~\ref{propro}, we obtain
the following proposition.

\begin{proposition}
$\displaystyle C_p(n,m)\ge \sum_{i=1}^{m}\sum_{a_1+\cdots+a_i=n-1}
\binom{n-1}{a_1,\ldots,a_i}.$
\end{proposition}

As a final remark, we observe, that another way to get the number
of different word patterns of length $n$ on $m$ letters is using a
correction in the formula for $T_p(m,n)$ like we did when we
obtained the upper bound for $C_p(n,m)$. But in this case the
correction is $i$ rather then $i-1$, which says that we consider
each conjugacy class with respect to shift and find the number of
elements in it. Thus, $i$ and $1/i$ cancel each other, and we get
a combinatorial proof of the following identity:
\[
\sum_{i=1}^{m}\sum_{a_1+\cdots+a_i=n \atop a_1\ge1,\dots,a_i\ge1}
\binom{n}{a_1,\ldots,a_i}= \sum_{i|n}
\sum_{j=0}^{\min(i,m)-1}(-1)^j\binom{\min(i,m)-1}{j}
\sum_{d|i}\mu(d)(\min(i,m)-j)^{\frac{i}{d}}.
\]


\end{document}